\documentclass[11pt]{article}

\usepackage{enumerate}
\usepackage{amsmath,xspace,amssymb,mathrsfs,theorem}

\input xy
\xyoption{all}
\xyoption{2cell}
\UseAllTwocells
\CompileMatrices

\title{Enriched weakness \thanks{Both authors gratefully acknowledge the support of the Ministry of Education of the Czech Republic by the project MSM 00216224409; the first-named author also gratefully acknowledges the support of the Australian 
Research Council.} }
\author{
Stephen Lack \\
Mathematics Department \\
Macquarie University NSW 2109 \\
Australia\\
{\tt steve.lack@mq.edu.au}
\and 
Ji\v r\'i Rosick\'y \\
Department of Mathematics and Statistics \\
Masaryk University
Kotl\'a\v rsk\'a 2 60000 Brno\\
Czech Republic\\
{\tt rosicky@math.muni.cz}
}
\date{}

\renewcommand{\phi}{\varphi}

\newcommand{\A}{{\ensuremath{\mathscr{A}}}\xspace}
\newcommand{\B}{{\ensuremath{\mathscr{B}}}\xspace}
\newcommand{\C}{{\ensuremath{\mathscr{C}}}\xspace}
\newcommand{\D}{{\ensuremath{\mathscr{D}}}\xspace}

\newcommand{\E}{{\ensuremath{\mathscr{E}}}\xspace}
\newcommand{\F}{{\ensuremath{\mathscr{F}}}\xspace}
\newcommand{\G}{{\ensuremath{\mathscr{G}}}\xspace}

\newcommand{\I}{{\ensuremath{\mathscr{I}}}\xspace}
\newcommand{\J}{{\ensuremath{\mathscr{J}}}\xspace}
\newcommand{\K}{{\ensuremath{\mathscr{K}}}\xspace}

\newcommand{\LL}{{\ensuremath{\mathscr{L}}}\xspace}

\newcommand{\T}{{\ensuremath{\mathscr{T}}}\xspace}

\newcommand{\V}{{\ensuremath{\mathscr{V}}}\xspace}

\newcommand{\X}{{\ensuremath{\mathscr{X}}}\xspace}

\newcommand{\Set}{\textnormal{\bf Set}\xspace}

\newcommand{\Cat}{\textnormal{\bf Cat}\xspace}

\newcommand{\nhom}{\textnormal{\bf NHom}\xspace}

\newcommand{\op}{\ensuremath{^{\textnormal{op}}}}
\newcommand{\fib}{\ensuremath{_{\textnormal{fib}}}}

\newcommand{\Pure}{\textnormal{Pure}}
\newcommand{\Ps}{\textnormal{Ps}}
\newcommand{\Cosk}{\textnormal{Cosk}}
\newcommand{\Iso}{\textnormal{Iso}}
\newcommand{\cod}{\textnormal{cod}}
\newcommand{\colim}{\textnormal{colim}}

\newcommand{\x}{\times}
\newcommand{\ot}{\otimes}

\newcommand{\ct}{\pitchfork}
\newcommand{\Inj}{\textnormal{Inj}}

\newdir{ >>}{{}*!/-10pt/@{>>}}
\newdir{ >}{{}*!/-8pt/@{>}}

\newcommand{\wrt}{with respect to }

\def\endproof{{\parfillskip=0pt\hfill$\Box$\vskip 10pt}}

\newcommand{\two}{\ensuremath{\mathbf{2}}\xspace}

\newtheorem{theorem}{Theorem}[section]    
\newtheorem{corollary}[theorem]{Corollary}   
\newtheorem{proposition}[theorem]{Proposition}   
\newtheorem{lemma}[theorem]{Lemma}   
\newtheorem{definition}[theorem]{Definition}

\newtheorem{schema}{Theorem-Schema}

{\theorembodyfont{\rmfamily}
\newtheorem{remark}[theorem]{Remark}   
\newtheorem{example}[theorem]{Example}   
}

\newcommand{\proof}{\noindent{\sc Proof:}\xspace}

\newcommand{\cotensor}{power\xspace}
\newcommand{\cotensors}{powers\xspace}

\begin{document}

\label{firstpage}
\maketitle

 \begin{abstract}
The basic notions of category theory, such as limit, adjunction, and orthogonality, 
all involve assertions of the existence and 
uniqueness of certain arrows. Weak notions arise when one drops the uniqueness 
requirement and asks only for existence. The enriched versions of the usual notions
involve certain  morphisms between hom-objects being invertible; here we introduce
enriched versions of the weak notions by asking that the morphisms between
hom-objects belong to a chosen class of ``surjections''.  We study in particular
injectivity (weak orthogonality) in the enriched context, and illustrate how it can be
used to describe homotopy coherent structures. 
 \end{abstract}

The basic notions of category theory, such as limit, colimit, free object, adjunction,
and factorization system, involve assertions of the existence of a unique morphism with certain properties. 
For example an object $FX$ is free on $X$ with respect to a functor $U:\A\to\X$
when there is a morphism $\eta:X\to UFX$, as in the diagram
$$\xymatrix{
X \ar[r]^{\eta} \ar[dr]_{f} & UFX \ar@{.>}[d]^{Ug} & FX \ar@{.>}[d]^{\exists! g} \\
& UA & A }$$
 with the property that 
for any morphism $f:X\to UA$ there is a unique morphism $g:FX\to A$ such
that $Ug.\eta=f$. It turns out that these can be expressed by saying that certain
induced functions between hom-sets are invertible --- in this case by saying that
the function $\A(FX,A)\to\X(X,UA)$ obtained by applying $U$ and then composing
with $\eta:X\to UFX$ is invertible for all $A\in\A$. It is this formulation in terms of hom-sets which makes 
enriched category theory possible --- you can replace
these bijections of hom-sets with isomorphisms of hom-objects lying in the
monoidal category \V over which categories are being enriched.

Weak notions arise when one asks just for the existence of a morphism with 
given properties, not the uniqueness. Thus in the example considered above,
we might ask that for every $f:X\to UA$ there exists a $g:FX\to A$ with
$Ug.\eta =f$. Once again this can be reformulated; this time by saying that the 
induced function $\A(FX,A)\to\X(X,UA)$ is {\em surjective} for all $A\in\A$.

In order to obtain an enriched version of these weak notions, it is necessary
to choose a class of morphisms in \V playing the role of the surjective 
functions. There are various possibilities: one might consider the epimorphisms,
or the regular epimorphisms, or the split epimorphisms; in fact we shall develop
the basic notions using an abstract class \E of morphisms in \V, although to make
much progress we shall have to start making some assumptions about this class.

This work began as a more-or-less technical investigation, within the context 
of our broader investigation of the homotopy theory of enriched categories. 
The role of weakness in homotopy theory is well-known: for instance 
weak factorization systems play a key role in the theory of model categories, 
and indeed weak limits were first considered in the homotopy context. In general,
however, weak limits are not the same as homotopy limits, and in homotopical 
situations, it is generally the latter that are more important. Nonetheless, as we
began to develop the theory and some specific examples, the project grew to be 
more than just technical. A key turning point was the example $\V=\Cat$ 
where \E consists of the {\em retract equivalences}. 

Since retract equivalences of categories are not just surjections
but rather some ``homotopy'' version of isomorphisms, this example
looks much less like ``weak category'' theory, as ordinarily understood, and 
much more like some sort of homotopy theory.  So in fact in this instance 
we are brought closer to a different, more recent, use of the word ``weak''
in category theory, meaning something like ``up to coherent homotopy''. We
illustrate this in Section~\ref{sect:Segal} by sketching how certain sorts of 
homotopy coherent structures can be described in terms of \E-injectivity classes for 
a suitable class \E.

The various weak notions we shall define will all say that some naturally
defined map or collection of maps lies in the class \E. Thus the smaller the 
class \E, the stronger the notion. In particular, the smallest possible class \E
is just the isomorphisms, and then our notion of weakness is not really weak at all,
and we recover the classical (non-weak!) theory of enriched categories.
The larger \E becomes, the weaker our weak notions really are.

The goal then, is to develop the basic ingredients of weak category theory,
such as can be found in \cite[Section~4A]{AR} for example, in this setting 
of \V-enriched category theory with a specified class \E of morphisms in \V
giving the weakness. In particular, we shall look at weak colimits, weak
adjunctions, injectivity,  and the basic
relationships between these notions. We hope to use this later in developing the 
homotopy theory of enriched categories, and in particular a homotopy version
 of locally presentable enriched categories.

The examples we have in mind are:

\begin{enumerate}[(a)]
\item If we take \E to be the isomorphisms, we obtain the ``non-weak
notion of weakness'': weak colimits are ordinary colimits, weak adjunctions are 
ordinary adjunctions, injectivity is orthogonality, and so on.
\item The ordinary notion of weakness, for $\V=\Set$, is where \E is the surjective
functions. We generalize this to the case of a locally
finitely presentable closed category \V by taking \E to be the pure
epimorphisms, whose definition is recalled below.
\item If $\V=\Cat$, we may take \E to be the retract equivalences.
\item If \V is a monoidal model category, we may take \E to be the trivial
fibrations.
\item \V is the category of fibrant objects in a monoidal model category, we may 
take \E to be  the weak equivalences.
\end{enumerate}

We prove our basic results for the first three cases. In the non-weak setting
of (a) all of this is known. In case (b) it is known for the case $\V=\Set$, but
new in general. Case (c) is new. Although we have  uniform statements
of results across these three cases, we have not been able to unify completely 
the proofs: while certain 
parts of the arguments are completely formal, others are treated on a case-by-case
basis. Case (d) is in fact a common generalization of the other three cases. We
are not able to prove all our results in the level of generality given in (d), but it would
be interesting to find conditions on a monoidal model category for which the
results can be proved.  (In fact we are not really using the full model structure, we
only use a weak factorization system: in this case the one involving the cofibrations
and trivial fibrations.)

 Case (e) is something like the setting of our broader (and on-going)
investigations into the homotopy theory of enriched categories. It is much more
delicate, since in general the fibrant objects in a monoidal model category will 
not be complete or cocomplete, or even closed under the monoidal structure. 
We shall have nothing to say about this case here, and indeed the approach
presented here has to be reformulated to deal with this case; nonetheless this case has 
influenced much of what we present below.

We start, in Section~\ref{sect:V}, by recalling a few key facts about enriched
category theory. We then make the basic definitions in the general setting:
injectivity in Section~\ref{sect:inj}, weak left adjoints in Section~\ref{sect:wla},
and weak colimits in Section~\ref{sect:wcolim}; in the latter we also define
a weakly locally presentable category to be an accessible category with 
weak colimits. These are the main objects of study. In the classical case, the
following conditions on a category \K are equivalent:
\begin{enumerate}[$(i)$]
\item \K is accessible and has weak colimits;
\item \K is accessible and has products;
\item \K is the full subcategory of a presheaf category consisting of 
those objects injective with respect to a given small set of morphisms;
\item \K is a weakly reflective subcategory of a presheaf category,
closed under retracts and under $\lambda$-filtered colimits for 
some regular cardinal $\lambda$;
\item \K is the category of models of a limit-epi sketch.
\end{enumerate}

We give analogous characterizations in the enriched context in each 
of the three main cases (a), (b), and (c) listed above.  This is done
in the three Sections~\ref{sect:isos}, \ref{sect:pure}, and \ref{sect:equivalence}.
Before this, in Section~\ref{sect:results}, we describe the general form of
our results, and in Section~\ref{sect:general} the common parts of the proofs.
We shall need a technical result involving enriched accessible categories; we 
state and prove it in Section~\ref{sect:technical}.
In Section~\ref{sect:Segal}, we describe various  examples of 
weakly locally presentable 2-categories, using $\V=\Cat$ and \E the retract
equivalences. These arise using the approach to coherent structures 
initiated by Segal in \cite{Segal74}. 
We describe in detail an example using bicategories.

\section{Our base category \V}\label{sect:V}

\enlargethispage{\baselineskip}
\enlargethispage{\baselineskip}

We work over a complete and cocomplete symmetric monoidal 
closed category \V, with tensor $\ot$ and unit $I$. We follow the general
conventions of \cite{Kelly-book}, and write $\V_0$ for the underlying ordinary category of \V. Later on we shall suppose that \V is locally finitely presentable as a closed category in the sense of Kelly \cite{Kelly-amiens}: this means that 
$\V_0$ is locally finitely presentable, the unit object $I$ is finitely presentable,
and the tensor product of two finitely presentable objects is finitely presentable.
It is not hard to generalize to any locally presentable $\V_0$.

We shall use heavily the notion of power (cotensor). For an object $A$ in
a \V-category \K and an object $X$ of \V, the power $X\ct A$ is defined
by the universal property
$$\K(-,X\ct A)\cong \V(X,\K(-,A)).$$
Dually the copower (or tensor) of an object $A\in\K$ by $X\in\V$  is
an object $X\cdot A$ defined by the universal property
$$\K(X\cdot A,-)\cong \V(X,\K(A,-)).$$

All notions should be understood as \V-enriched unless specified otherwise.
For instance,  if we speak of a full subcategory or a functor, this should always be
understood to mean a full sub-\V-category or \V-functor, as the case may be.
Similarly limits or colimits will always be understood in the enriched sense,
defined in terms of some isomorphism of \V-valued homs. 

If we do wish to speak of ordinary unenriched notions, we shall say so. We
shall sometimes consider ordinary diagrams in enriched categories: given
a \V-category \K and an ordinary category \J, a diagram in \K of shape \J
is an ordinary functor from \J to the underlying ordinary 
category of \K (although this is equivalent to giving a \V-functor from the
free \V-category on \J to \K). The (conical) colimit of such a diagram $S$ in a \V-category
\K is defined by a natural isomorphism 
$$\K(\colim S,A)\cong [\J,\K](S,\Delta A)$$
in \V. As explained in \cite{Kelly-book}, the universal property of a conical colimit
in the underlying ordinary category $\K_0$ is weaker than  that of a colimit in \K;
but if the colimit in \K is known to exist, then the universal property of the colimit in 
$\K_0$ is enough to detect it. As a consequence, if \K and \LL have conical colimits
of a certain type, then a \V-functor $F:\K\to\LL$ will preserve these colimits 
provided that the underlying ordinary functor $F_0:\K_0\to\LL_0$ does so.

 Furthermore, there is an important special case in which
the universal property of the colimit in $\K_0$ suffices: let \G be a strong generator for $\V_0$, and suppose that \K has powers by all $G\in\G$. Then there is no difference
between the colimit of $S$ in \K, and the colimit in $\K_0$. In particular,
this is the case if $\V_0$ is locally finitely presentable and \K has powers by
finitely presentable objects of $\V_0$.

We follow \cite{Kelly-book} by using {\em filtered colimit} to mean 
the (conical) colimit of a diagram $\J\to\K$ with \J an ordinary filtered category;
the case of $\lambda$-filtered colimits, for a regular cardinal $\lambda$, is similar.
We say that a \V-category \A is {\em $\lambda$-accessible}, if it is the free
completion under $\lambda$-filtered colimits of a small \V-category \C. 

Since conical colimits can be defined in terms of the underlying categories
in the presence of enough powers, we have:

\begin{proposition}
If a \V-category  has all \G-powers
for some strong generator \G of $\V_0$, then \A is
$\lambda$-accessible if and only if the underlying ordinary category $\A_0$
is so. 
\end{proposition}

\begin{remark}
The words {\em filtered} and {\em accessible} were used in a different sense,
in relation to \V-categories, in \cite{BQR:V-accessible}; there filtered colimits
included all weighted colimits commuting with finite limits, and accessibility
was defined accordingly.
\end{remark}

We shall also fix a class \E of morphisms in \V, which  to start with is 
assumed only to be closed under composition and to contain the isomorphisms.
As explained in the introduction, ``weak category theory'' corresponds to 
the case $\V=\Set$ and \E the surjections, while the non-weak case corresponds
to taking \E to consist of  the isomorphisms.

\section{Injectivity in \V-categories}\label{sect:inj}

\enlargethispage{\baselineskip}

Let $f:A\to B$ be a morphism in a \V-category \K. We say that an object
$C\in\K$ is {\em \E-injective with respect to $f$}, or just {\em injective with respect to $f$} if \E is understood,
when the induced morphism 
$$\xymatrix{ \K(B,C) \ar[r]^{\K(f,C)} & \K(A,C) }$$
is in \E. More generally, if \F is any class of morphisms in \K, we say that 
$C$ is \E-injective with respect to \F  if it is \E-injective with respect to  
all $f\in\F$. We write
$\Inj_\E(\F)$ for the full subcategory of \K consisting of those objects which 
are \E-injective with respect to  all $f\in\F$, and call such a full subcategory 
an {\em injectivity class}, or {\em \E-injectivity class} for emphasis. We 
call it a small-injectivity class if the class \F of morphisms is small.
We write $\Inj_0(\F)$
for the full subcategory of \K consisting of those objects which are injective
with respect to $f$ in the ordinary category $\K_0$ in the ordinary sense, 
for all $f\in\F$.

Obviously this notion of injectivity depends heavily on both \V and \E,
as the following examples show:

\begin{example}
Ordinary injectivity is the case $\V=\Set$ and \E the surjections.  
\end{example}

\begin{example}
We can also obtain {\em orthogonality} as an example, by taking \E to 
be the isomorphisms. This works in either the enriched or the unenriched
contexts.  
\end{example}

\begin{example}
Let $\V=\Cat$ and \E be the equivalences. Let $1$ and $2$ be discrete categories
with one and two objects, respectively, and let \I be the free-living isomorphism,
consisting of two objects and an invertible arrow between them. Let $f:1\to 2$ be
an injection, and $g:\I\to 1$ the unique functor into the terminal category. Since
$g$ is an equivalence, all objects of \Cat are \E-injective with respect to $g$; 
on the other 
hand very few categories are injective in $\Cat_0$ with respect to $g$ --- in 
particular, \I is not.
Since $f$ is split mono, all categories are injective in $\Cat_0$ with respect to $f$, 
while very few are \E-injective with respect to $f$ --- in particular, $2$ is not. 
\end{example}

We do have the following positive result, which requires the unit object $I$ of \V
to be {\em \E-projective}, in the sense that $\V_0(I,-):\V_0\to\Set$ sends morphisms
in \E to surjections; more explicitly, this says that for any $e:X\to Y$ in \E and 
any $y:I\to Y$, there exists an $x:I\to X$ such that $ex=y$. 

\begin{proposition}
If $I$ is \E-projective, then \E-injectivity implies ordinary injectivity in the 
underlying ordinary category.
\end{proposition}

\proof
To say that $C$ is \E-injective with respect to $f$ is to say that the morphism 
$\K(f,C):\K(B,C)\to\K(A,C)$ in \V lies in \E. Under the conditions of the 
proposition, however, this implies that $\V_0(I,\K(f,C))$ is surjective.
Now $\V_0(I,\K(f,C))$ is just 
$\K_0(f,C):\K_0(B,C)\to\K_0(A,C)$, and surjectivity of 
$\K_0(f,C)$ is just ordinary injectivity of $C$ with respect to $f$.
\endproof

\noindent
This condition on the unit object will hold in all the main examples we study.
We shall see in Proposition~\ref{prop:Inj0} that in many cases \E-injectivity 
with respect to \F is 
equivalent to 
ordinary injectivity with respect to some other class $\F'$ of maps. In this 
case, if also $I$ is \E-projective, then  
\E-injectivity classes reduce to injectivity classes in the usual sense. 

\begin{definition}\label{defn:E-stable}
We say that a class of limits is {\em \E-stable} if \E is closed in 
$\V^\two$ under these limits.
\end{definition}

\begin{proposition}\label{prop:lims}
Any \E-injectivity class $\Inj_\E(\F)$ in \K 
 is closed under \E-stable limits.
\end{proposition}

\proof
Let $S:\D\to\K$ be any diagram in \K for which $SD$ is injective with respect to $f$
for all $D\in\D$ and all $f:A\to B$ in \F. 
We consider a limit (possibly weighted) of $S$. In the
diagram
$$\xymatrix{
\K(B,\lim S) \ar[r] \ar[d]_{\K(f,\lim S)} & \lim \K(B,SD) \ar[d]^{\lim \K(f,SD)} \\
\K(A,\lim S) \ar[r] & \lim \K(A,SD) }$$
the horizontal arrows are invertible, since representables preserve
limits. Thus the left hand vertical is in \E if the right hand vertical is.
\endproof

\begin{proposition}\label{prop:colims}
The full subcategory $\Inj_\E(\F)$ of \K consisting of the objects injective with
respect to \F is also closed under any class $\Phi$  of colimits for which 
\begin{enumerate}[$(i)$]
\item \E is closed under $\Phi$-colimits
\item $\K(A,-)$ preserves $\Phi$-colimits for any object $A$ which is 
the domain or codomain of a morphism in \F.
\end{enumerate}
\end{proposition}

\proof
Let $S:\D\to\K$ be any diagram in \K for which $SD$ is injective with respect to $f$
for all $D\in\D$ and all  $f:A\to B$ in \F. 
We consider a limit (possibly weighted) of $S$. In the
diagram
$$\xymatrix{
\colim \K(B,SD) \ar[d]_{\colim \K(f,SD)} \ar[r] & \K(B,\colim S)  \ar[d]^{\K(f,\colim S)}  \\
\colim \K(A,SD) \ar[r] & \K(A,\colim S)  }$$
the horizontal arrows are invertible, since $\K(B,-)$ and $\K(A,-)$ are assumed
to preserve the colimits in question, and the left hand vertical is in \E,
since the \E's assumed to be closed under the colimits in question. Thus
the right hand vertical is in \E.
\endproof

\begin{remark}
It is not hard to generalize \E-injectivity to a notion of \E-weak
factorization system. 
\end{remark}

\section{Weak left adjoints}\label{sect:wla}

\enlargethispage{\baselineskip}
\enlargethispage{\baselineskip}

Let $U:\A\to\K$ be a \V-functor, and $K$ an object of \K. We
say that a morphism $\eta:K\to UFK$ exhibits $FK$ as a weak
left adjoint to $U$ at $K$, when
for any $A\in\A$ the induced map
$$\xymatrix @C1.5pc {
\A(FK,A) \ar[r]^-{U} & \K(UFK,UA) \ar[rr]^-{\K(\eta,UA)} && \K(K,UA) }$$
is in \E.

A special case is where we actually have a functor $F:\K\to\A$, and
a natural transformation $\eta:1\to UF$, 
for which 
$$\xymatrix @C1.5pc {
\A(FK,A) \ar[r]^-{U} & \K(UFK,UA) \ar[rr]^-{\K(\eta,UA)} && \K(K,UA) }$$
is in \E.
We then say that $U$ has a {\em natural weak left adjoint}. 


If there is a weak left adjoint to $U$ at every object of \K, we say 
simply that $U$ has a weak left adjoint. We shall be particularly interested in the 
case where $U$ is fully faithful and has a weak left adjoint, in which case we say 
that \A is {\em weakly reflective} in \K (or \E-weakly reflective if we wish to 
emphasize \E).

\begin{example}
Suppose that \E contains the retractions.
If $W:\LL\to\K$ has a left adjoint $H$, and we factorize $W$ as a 
bijective on objects functor $P:\LL\to\A$ followed by a fully
faithful $U:\A\to\K$, then let $F=PH$ and $\eta:1\to UF=WH$
be the unit of the adjunction. We claim that this exhibits $F$
as a natural weak left adjoint (a natural weak reflection) to $U$.
To see this, we shall show that 
$$\xymatrix @C1.5pc {
\A(FK,A) \ar[r]^-{U} & \K(UFK,UA) \ar[rr]^-{\K(\eta,UA)} && \K(K,UA) }$$
has a section. To do this, observe that each $A\in\A$ has
the form $PL$ for a unique $L\in\LL$, and that the counit 
$\epsilon L:HWL\to L$ gives a map $P\epsilon L:HUPL=HWL\to PL$
so that we obtain (non-natural) maps $\pi A:HUA\to A$ 
with $\pi PL=P\epsilon L$.  Now the required section is given by
$$\xymatrix @C1.5pc {
\K(K,UA) \ar[r]^-{F} & \A(FK,FUA) \ar[rr]^-{\A(FK,\pi A)} && \A(FK,A) }.$$
\end{example}

The following proposition says roughly that weak left adjoints compose.

\begin{proposition}
Let $U:\A\to\B$ and $V:\B\to\C$ be \V-functors. Suppose that 
$\eta:C\to VGC$ exhibits $GC$ as a weak left adjoint to $V$ at $C$,
and $\beta:GC\to UFGC$ exhibits $FGC$ as a weak left adjoint to $U$
at $GC$. Then the composite
$$\xymatrix{
C \ar[r]^-{\eta} & VGC \ar[r]^{V\beta} & VUFGC }$$
exhibits $FGC$ as a weak left adjoint to $VU$ at $C$.
\end{proposition}

\proof
For any $A\in\A$ we have 
$$\xymatrix{
\A(FGC,A) \ar[r]^{U} \ar[dr] & 
\B(UFGC,UA) \ar[r]^{V} \ar[d]^{\B(\beta,UA)} &
\C(VUFGC,VUA) \ar[d]^{\C(V\beta,VUA)}\\
& \B(GC,UA) \ar[r]^{V} \ar[dr] & \C(VGC,VUA) \ar[d]^{\C(\eta,1)} \\
&& \C(C,VUA) }$$
and the two diagonals are in \E so the composite diagonal is in \E.
\endproof

There is a corresponding result for natural weak left adjoints.

\enlargethispage{\baselineskip}

\section{Weak colimits}\label{sect:wcolim}

\enlargethispage{\baselineskip}
\enlargethispage{\baselineskip}

Let $S:\C\to\K$ be a \V-functor. There is an induced \V-functor
$$\xymatrix{
\K \ar[r]^-{\K(S,1)} & [\C\op,\V] }$$
sending an object $A\in\K$ to $\K(S-,A)$. We sometimes write
$\widetilde{S}$ for $\K(S,1)$.

Now let $H:\C\op\to\V$ be a \V-functor.
We define a {\em weak $H$-weighted colimit} of $S$ to be a weak left adjoint 
to $\K(S,1)$ at $H$. Explicitly, this means we have an object $C\in\K$
and a natural transformation $\gamma:H\to\K(S,C)$, such that for all $A\in\K$
the induced map 
$$\xymatrix{
\K(C,A) \ar[r]^-{\K(S,1)} & [\C\op,\V](\K(S,C),\K(S,A)) \ar[rr]^{ [\C\op,\V](\gamma,1)} &&
[\C\op,\V](H,\K(S,A)) }$$
in \V is in \E. We may sometimes write $H*_\E S$ or $H*_w S$ for such a weak 
colimit.

\begin{example}
If \E is the class of isomorphisms, then this reduces to the usual notion of weighted 
colimit.
\end{example}

\begin{example}
If $\V=\Set$ and \E is the class of surjections, a weak colimit $H*_w S$ consists of 
an object $C$ equipped with a natural transformation $\gamma:H\to\K(S,C)$
such that for any $A$ and any natural transformation $\alpha:H\to\K(S,A)$ 
there exists a morphism $f:C\to A$, not necessarily unique, such that
$\K(S,f)\gamma=\alpha$. In particular, this reduces to ordinary weak (conical)
colimits when the weight $H:\C\op\to\Set$ is constant at the terminal object 
of \Set.
\end{example}

Recall that $S:\C\to\K$ is said to be dense when the induced map 
$\K(S,1):\K\to[\C\op,\V]$, sending $A\in\K$ to $\K(S-,A):\C\op\to\V$,
is fully faithful.

\begin{proposition}\label{prop:weak-cocompleteness}
If $S:\C\to\K$ is dense and \C small, then 
$\K(S,1)$ has a weak left adjoint if and only if \K has weak colimits.
\end{proposition}

\proof
By definition, $\K(S,1)$ has a weak left adjoint if \K has all weak
colimits $H*_w S$. This gives one direction. For the converse, suppose that
$\K(S,1)$ has a weak left adjoint, and that $R:\A\to\K$ is any (small)
diagram. Form the composite $\K(S,R):\A\to[\C\op,\V]$ and consider
the composite 
$$\xymatrix{
\K \ar[r]^-{\K(S,1)} & [\C\op,\V] \ar[r]^-{\widetilde{\K(S,R)}} & [\A\op,\V]}$$
where the second functor sends $G:\C\op\to\V$ to the functor
$\A\op\to\V$ sending $A\in\A$ to $[\C\op,\V](\K(S-,RA),G-)$.
Now $\K(S,1)$
has a weak left adjoint by assumption, while $\widetilde{\K(S,R)}$
has an actual left adjoint, since $[\C\op,\V]$ is cocomplete. So the composite has a 
weak left adjoint. But
$\K(S,1)$ sends an object $A\in\K$ to $\K(S,A)$, and $\widetilde{\K(S,R)}$ now sends
this $\K(S,A)$ to $[\C\op,\V](\K(S,R),\K(S,A))$; but since $\K(S,1)$ is fully faithful,
this is just $\K(R,A)$, and so the composite is really just $\K(R,1)$. 
We have therefore shown that $\K(R,1)$ has a weak left adjoint, and 
so that \K has weak colimit $H*_w R$ for all weights $H$.
\endproof

As in the case of \Set, we define a \V-category \A to be {\em weakly 
locally presentable} if it is accessible and has weak colimits.

 In the 
case of \Set, the weakly locally presentable categories are precisely 
the categories of models of limit-epi sketches; we shall also prove
enriched versions of this characterization. For convenience, the 
limit part of our sketches will be taken to be limit theories (that is,
small \V-categories with $\alpha$-small limits for some regular
cardinal $\alpha$).

A {\em (limit,\E)-sketch} is a small \V-category \T with $\alpha$-small
limits for some regular cardinal $\alpha$, and a specified collection \F of 
morphisms in \T. The \V-category
of models of the sketch is the full subcategory of $[\T,\V]$ consisting
of those \V-functors which preserve $\alpha$-small limits and send
morphisms in \F to morphisms in \E. In practice, the category \T 
with $\alpha$-small limits may be presented via a sketch.

\section{The structure of the theorems}\label{sect:results}

\enlargethispage{\baselineskip}
\enlargethispage{\baselineskip}

In this section we explain briefly the form of our results. At this stage 
we suppose only that the monoidal category \V is locally presentable, 
which we henceforth assume. We describe the results in the form of three 
Theorem-Schemas, but should point out straight away that they do not hold
without further assumptions on \V and \E. 
In this section we show that the second and third of these Theorem-Schema
follow from the first. 
In the sections that follow, we shall describe the various examples in which
we can prove the first.

\begin{schema}\label{schema:main}
Let \K be a locally presentable \V-category, and \A a full subcategory.
The following conditions are equivalent:
\begin{enumerate}[$(i)$]
\item \A is the category of objects \E-injective to a small class of maps in \K (a 
small-\E-injectivity class in \K);
\item \A is accessible, accessibly embedded, and closed under \E-stable limits;
\item \A is accessibly embedded and \E-weakly reflective.
\end{enumerate}
\end{schema}

Following \cite{AR}, we say that a full subcategory of an accessible category 
is {\em accessibly embedded} if it is closed under $\alpha$-filtered colimits for 
some regular cardinal $\alpha$. We understand this to include the fact that
the subcategory is {\em replete}, meaning that it if it contains an object $A$ then it 
contains any object isomorphic to $A$.

\begin{schema}
For any \V-category \A, the following conditions are equivalent:
\begin{enumerate}[$(i)$]
\item \A is an \E-weakly reflective, accessibly embedded, full subcategory
of some presheaf category $[\C,\V]$;
\item \A is (equivalent to) 
a small-\E-injectivity class in some locally presentable \V-category \K;
\item \A is accessible and has \E-stable limits
\item \A  is accessible and has \E-weak colimits.
\end{enumerate}
\end{schema}
A \V-category \A satisfying these conditions is said to be 
{\em \E-weakly locally presentable}, or just {\em weakly 
locally presentable} if \E is understood.

\proof 
We shall prove the equivalence, given that Theorem-Schema~\ref{schema:main} holds.
The implications $(i)\Rightarrow(ii)\Rightarrow(iii)$ follow immediately 
from Theorem-Schema~\ref{schema:main}. 
To see that $(iii)\Rightarrow(i)$, suppose that \A
is accessible and has \E-stable limits. Choose a regular cardinal $\lambda$
for which \A is $\lambda$-accessible, and let $\A_\lambda$ be the full
subcategory of \A consisting of the $\lambda$-presentable objects. Then \A
is a full subcategory of $[\A\op_\lambda,\V]$, closed under $\lambda$-filtered
colimits and \E-stable limits. By Theorem-Schema~\ref{schema:main} it is \E-weakly
reflective, and so $(i)$ holds.

If the first three conditions hold, then since \A is \E-weakly reflective in 
the cocomplete $[\C,\V]$, it is \E-weakly cocomplete. Thus $(iv)$ holds.

Conversely, if \A is accessible with \E-weak colimits, choose a regular cardinal
$\lambda$ such that \A is $\lambda$-accessible, and consider the embedding
$\A\to[\A\op_\lambda,\V]$. Since \A is \E-weakly cocomplete, it is \E-weakly
reflective by Proposition~\ref{prop:weak-cocompleteness}. This gives $(i)$.
\endproof

Finally there is a description in terms of sketches. We shall not consider
the most general notion of sketch, but restrict ourselves to the case of a 
small \V-category \T with certain specified limit
diagrams, and with a chosen class $\F$ of morphisms. A model is
then a \V-functor from \T to \V which sends the specified limits to limits in \V, and 
the morphisms in $\F$ to morphisms in \E.

\begin{schema}
A \V-category \A is \E-weakly locally presentable if and only if it is equivalent to 
the category of models of a (limit,\E)-sketch.
\end{schema}

\proof
Once again, we prove this given that Theorem-Schema~\ref{schema:main} holds. 
Suppose first that \T is a small \V-category with certain specified limit diagrams,
and that \F is a set of morphisms in \T. Let $\K$ be the full subcategory of $[\T,\V]$
sending the chosen limits in \T to limits in \V. Then \K is locally presentable. 
Now a \V-functor $M:\T\to\V$ in \K sends  the morphisms in \F to \E
if and only if it is injective in \K with respect to the morphisms $\T(f,-):\T(B,-)\to\T(A,-)$ for all $f:A\to B$ in \F. Thus the models of a (limit,\E)-sketch are a small-injectivity
class.

For the converse, suppose that $\A=\Inj_\E(\F)$ for some small set \F of morphisms
in a locally presentable \V-category \K. Choose a regular cardinal $\lambda$ 
sufficiently large that \K is locally $\lambda$-presentable, and all domains and
codomains of morphisms in \F are $\lambda$-presentable in \K. Let \C be 
the opposite of the category of $\lambda$-presentable objects in \K. Then \C
has $\lambda$-small limits, and \K is equivalent to the category of $\lambda$-continuous functors from \C to \V. Furthermore, \F can be seen as a set of morphisms in 
\C, and the objects of \A are those objects of \K which are \E-injective with respect
to the morphisms in \F. Thus $(\C,\F)$ is a (limit,\E)-theory whose \V-category of
models is \A.
\endproof

\section{General aspects of the proofs}\label{sect:general}

\enlargethispage{\baselineskip}

In this section we introduce further assumptions which allow us to prove the 
equivalences in Theorem-Schema~\ref{schema:main}. The first assumption is 
easy: we suppose that $I$ is \E-projective; recall that this means that 
$\V_0(I,-):\V_0\to\Set$ sends the maps in \E to surjections.

\begin{proposition}
If $I$ is \E-projective then the
implication $(iii)\Rightarrow(i)$ in Theorem-Schema~\ref{schema:main} holds.
\end{proposition}

\proof
For each object $K\in\K$, we may choose a weak reflection $r_K:K\to K^*$ into \A.
The universal property of the weak reflection says that each object $A\in\A$ is \E-injective with respect
to these maps $r_K$. Conversely, if $K\in\K$ is \E-injective with respect to the
single map $r_K:K\to K^*$, then since $I$ is projective with respect to \E, 
it follows that $K$ is injective in $\K_0$ with respect to $r_K$, and so that $K$
is a retract of $K^*$. But \A is accessibly embedded, so closed under
retracts, thus $K\in\A$.

Thus \A consists of all objects which are \E-injective with respect to all the $r_K$. 
The only problem is that this is a large class of maps; to prove the proposition,
we must show that it can be replaced by a small one. 

Choose a regular cardinal $\lambda$ such that \A is $\lambda$-accessible,
and the inclusion $\A\to\K$ preserves $\lambda$-filtered colimits and 
$\lambda$-presentable objects.
Let \F consist of all the $r_K:K\to K^*$ for which $K$ is $\lambda$-presentable
in \K. Certainly \A is contained in $\Inj_\E(\F)$; we must show that the reverse
inclusion holds.

Suppose then that $X\in\Inj_\E(\F)$. Let \J be the full subcategory of $\K_0/X$ 
consisting of all morphisms into $X$ with $\lambda$-presentable domain. Then
\J is $\lambda$-filtered, and $X$ is the colimit of the canonical map $\J\to\K$. 

Let $\J'$ be the full subcategory of $\J$ consisting of those $K\to X$ for which $K$
is not just $\lambda$-presentable, but also in \A. We shall show that $\J'$ is 
final in \J. Then $\J'$ will still be 
$\lambda$-filtered, and $X$ will be a $\lambda$-filtered colimit of objects
in \A, and so itself will be in \A. 

Since \J is $\lambda$-filtered, $\J'$ will be final provided that for each $J\in\J$,
there is a morphism $J\to J'$ with $J'\in\J'$. So let $f:K\to X$ in \J be given. 
Since \A is $\lambda$-accessible, $K^*$ is a $\lambda$-filtered colimit of
$\lambda$-presentable objects of \A. Since $K$ is $\lambda$-presentable,
$r:K\to K^*$ factorizes as $s:K\to B$ followed by $g:B\to K^*$ for some
$\lambda$-presentable object $B\in\A$. Since $r:K\to K^*$ is a weak
reflection, $f=f'r$ for some $f':K^*\to X$, and $f=f'r=f'gs$. We have an
object $f'g:B\to X$ of $\J'$, and a morphism $s$ from $f:K\to X$ to $f'g:B\to X$ in \J.
This proves that $\J'$ is final in \J, and so completes the proof.
\endproof

Next we introduce a condition that allows enriched injectivity to be reduced to
ordinary injectivity. We say that the class \E is {\em cofibrantly generated} if 
there is a small set \J of  morphisms in \V such that \E consists of those 
morphisms with the right lifting property with respect to \J; in other words, those
$e$ with the property that for any commutative square as in the solid part of the
diagram
$$\xymatrix{
\ar[d]_{j} \ar[r]^-u & \ar[d]^e \\
\ar[r]_-v \ar@{.>}[ur]_w & {} }$$
with $j$ in \J, there exists a ``fill-in'' $w$ making the two triangles commute. 
(It then follows that \E forms part of a cofibrantly generated weak factorization system.)

\enlargethispage{\baselineskip}
\enlargethispage{\baselineskip}

\begin{proposition}\label{prop:Inj0}
If \E is cofibrantly generated and \F is a small set of morphisms in \K then there
is a small set $\F'$ of morphisms in \K for which \E-injectivity with respect to  \F 
is equivalent to ordinary injectivity with respect to $\F'$; that is, $\Inj_\E(\F)=\Inj_0(\F')$.
\end{proposition}

\proof
Let \J be a set of morphisms which generates \E in the sense of the previous paragraph. 
An object $A\in\K$ is \E-injective \wrt \F if and only if $\K(f,A)$ is in \E for all
$f:B\to C$ in \F. But this says that $\K(f,A)$ has the right lifting property with 
respect to all $j$ in \J; in other words, given the solid part of the diagram 
$$\xymatrix{
X \ar[r]^-u \ar[d]_j & \K(C,A) \ar[d]^{\K(f,A)} \\
Y \ar[r]_-v \ar@{.>}[ur] & \K(B,A) }$$
there exists a diagonal fill-in. But to give $u$ and $v$ is equivalently to give
$u':X\cdot C\to A$ and $v':Y\cdot B\to A$ making the square 
$$\xymatrix{
X\cdot B \ar[r]^{X\cdot f} \ar[d]_{j\cdot B} & X\cdot C \ar[d]^{u'} \\
Y\cdot B \ar[r]_{v'} & A }$$
commute; or equivalently a morphism 
$$\xymatrix{
X\cdot C +_{X\cdot B} Y\cdot B \ar[r]^-{w} & A }$$
out of the pushout of $X\cdot f$ and $j\cdot B$. A diagonal fill-in is then 
equivalent to the existence of an extension of $w$ along the canonical map
$$\xymatrix{ X\cdot C+_{X\cdot B} Y\cdot B \ar[r]^-{j\Box f} & Y\cdot C.}$$
Thus $A$ is \E-injective \wrt \F if and only if it is injective \wrt $\F'$ in the ordinary
sense, where $\F'$ consists of all $j\Box f$ with $j\in\J$ and $f\in\F$.
\endproof

We now state a second implication from Theorem-Schema~\ref{schema:main};
the proof depends on the notion of pure subobject, but we defer this aspect to 
Section~\ref{sect:technical}.

\begin{corollary}\label{cor:one-two}
If \E is cofibrantly generated, then the implication $(i)\Rightarrow(ii)$ in 
Theorem-Schema~\ref{schema:main} holds.
\end{corollary}

\proof
Let \K be locally presentable, and \F a small set of morphisms in \K.
We have already seen that $\Inj_\E(\F)$ is closed under \E-stable limits and
under $\lambda$-filtered colimits, where $\lambda$ is a regular cardinal
for which the domain and codomain of any morphism in \F is $\lambda$-presentable.
It remains only to prove that $\Inj_\E(\F)$ is accessible.
Now the underlying ordinary category $\Inj_\E(\F)_0$ of $\Inj_\E(\F)$ is
$\Inj_0(\F')$ by the previous proposition, and by the classical theory this is an accessible
and accessibly embedded subcategory of the locally finitely presentable ordinary
category $\K_0$. Thus  by \cite[Corollary~2.36]{AR},  $\Inj_\E(\F)_0$ is closed in $\K_0$
under $\mu$-pure subobjects for some regular cardinal $\mu$, and now $\Inj_\E(\F)$
is accessible by Theorem~\ref{thm:purity}.
\endproof

Thus if \E is cofibrantly generated and $I$ is \E-projective then we have
the implications $(iii)\Rightarrow(i)\Rightarrow(ii)$, and it remains only to 
prove $(ii)\Rightarrow(iii)$. Rather than describing general sufficient conditions 
for this to be true, we treat it on a case-by-case basis in the examples that follow.

\section{Case 1: \E is the isomorphisms}\label{sect:isos}

This case is entirely classical, we merely state the results, to show what 
they give in this context. In this case injectivity becomes orthogonality,
weak reflectivity is ordinary reflectivity, all limits are \E-stable,
and weak colimits are just colimits.

Observe that $I$ is indeed \E-projective, since $\V_0(I,-)$ sends the isomorphisms
not just to surjections but to bijections (as indeed does any functor). Also the
class \E is cofibrantly generated: if $\alpha$ is some regular cardinal for which
$\V_0$ is locally $\alpha$-presentable, then the $\alpha$-presentable objects
form a strong generator for $\V_0$, and  a morphism $e:A\to B$ in $\V_0$ is in \E
(that is, invertible) if and only 
if it has the right lifting property with respect to the unique map $0\to G$ 
and the codiagonal $G+G\to G$ for all $\alpha$-presentable objects $G$.

\begin{theorem}
Let \K be a locally presentable \V-category, and \A a full subcategory. 
The following conditions are equivalent:
\begin{enumerate}[$(i)$]
\item \A is the category of objects orthogonal to a small class of maps in \K;
\item \A is accessible, accessibly embedded, and closed under limits;
\item \A is  accessibly embedded and reflective.
\end{enumerate}
\end{theorem}

By the general results of the previous section, we need only prove that (ii) implies
(iii). Suppose then that \A is accessible, accessibly embedded, and closed under
limits. Let $K\in\K$ be given; we shall construct a reflection of $K$ into \A. Let
$J:\A\to\K$ be the inclusion, and
$\K(K,J):\A\to\V$ be the \V-functor sending $A\in\A$ to the hom $\K(K,JA)$.
If $\lambda$ is any regular cardinal for which \A is closed in \K under $\lambda$-filtered colimits and $K$ is $\lambda$-presentable, then $\K(K,J)$ will preserve
$\lambda$-filtered colimits. We may choose $\lambda$ so that \A is also 
$\lambda$-accessible, and now $\K(K,J)$ is the left Kan extension of its restriction
to the $\lambda$-presentable objects. Thus $\K(K,J)$-weighted limits exist in 
\A and are preserved by $J$. The $\K(K,J)$-weighted limit of the identity 
$\A\to\A$ is now the desired reflection of $K$ into \A. 

The other two theorems now follow as in the previous section:

\begin{theorem}
For a \V-category \A, the following are equivalent:
\begin{enumerate}[$(i)$]
\item \A is a reflective, accessibly embedded subcategory of $[\C,\V]$ for some small \V-category \C;
\item \A is equivalent to a small orthogonality class in some locally presentable \V-category;
\item \A is accessible and complete;
\item \A is accessible and cocomplete.
\end{enumerate}
\A is then said to be locally presentable.
\end{theorem}

\begin{theorem}
A \V-category is  locally presentable  if and only if it is equivalent to the
\V-category of models of a limit  sketch.
\end{theorem}

The general form of this last theorem would be that a \V-category
is \E-weakly locally presentable if and only if it is the \V-category of
models of a (limit,\E)-sketch. But when \E is the isomorphisms \E-weakly means
not weakly at all; and \E-specifications are just iso-specifications, 
which do not require any sort of colimits.

\section{Case 2: \E is the pure epimorphisms}\label{sect:pure}

\enlargethispage{\baselineskip}
\enlargethispage{\baselineskip}

In this section, we suppose that \V is locally finitely presentable
as a closed category \cite{Kelly-amiens}; recall that this means that the underlying 
ordinary category $\V_0$ is locally finitely presentable, and the full subcategory 
of finitely presentable objects contains the unit and is closed under the tensor 
product. 

Recall that an epimorphism $p:X\to Y$ in a locally finitely presentable
(ordinary) category \K 
is said to be {\em pure} \cite{AR:pure}, if $\K(G,p):\K(G,X)\to\K(G,Y)$ is surjective
for all finitely presentable objects $G$.
We now take as our class \E the pure epimorphisms in $\V_0$: equivalently these are the morphisms with the right lifting property with respect to the 
unique map $0\to G$ for all finitely presentable objects $G$, so \E is cofibrantly 
generated. Note also that $I$ is finitely presentable, and so is \E-projective.
Thus all our theorems will hold provided that the implication 
$(ii)\Rightarrow(iii)$ in Theorem-Schema~\ref{schema:main} does, which we shall
see below.

As usual, we regard \E as a full subcategory of $\V^\two$.

\begin{proposition}
The pure epimorphisms are closed in \textrm{$\V^\two$} under products, retracts,
and finite powers.
\end{proposition}

\proof
Let $\Pi_i p_i:\Pi_i X_i\to \Pi_i Y_i$ be a product of pure epimorphisms,
let $G$ be finitely presentable, and let $f:G\to \Pi_i Y_i$ be given. Then 
$f$ is determined by components $f_i:G\to Y_i$ for each $i\in I$. Since
$p_i$ is pure epi and $G$ is finitely presentable, there is a lifting 
$f_i = p_i g_i$ of $f_i$ through $p_i$ for each $i$, and so a lifting 
of $f$ through $\Pi_i p_i$. This proves that the pure epimorphisms are
closed under products. The case of retracts is similar.

To see that the pure epimorphisms are closed under finite powers, let
$p:X\to Y$ be a pure epimorphism, and $H\in\V$ a finitely presentable 
object. We must show that $H\ct p:H\ct X\to H\ct Y$ is a pure epimorphism.
Suppose then that $f:G\to H\ct Y$ is given; this amounts to giving $f':G\ot H\to Y$.
Since $G$ and $H$ are finitely presentable, and \V is locally finitely presentable
as a closed category, $G\ot H$ is also finitely presentable, and so $f'$ lifts
through $p$, say as $g':G\ot H\to X$. Now $g'$ determines a unique 
$g:G\to H\ct X$, and $(H\ct p)g=f$, which proves that $H\ct p$ is a pure
epimorphism.
\endproof

\begin{remark}
There are analogues to all results in this section for higher
cardinals $\lambda$. This would involve a \V which is locally
$\lambda$-presentable as a closed category, and taking \E to be 
the $\lambda$-pure epis (that is, the morphisms $p:X\to Y$ with 
$\K(G,p)$ surjective for all $\lambda$-presentable $G$). 
We shall not bother to spell these out.
\end{remark}

\begin{proposition}[\cite{AR:pure}]
Pure epimorphisms are closed under filtered colimits.
\end{proposition}

\proof
Consider a filtered colimit $\colim_i p_i:\colim_i X_i\to \colim_i Y_i$ in 
$\V^\two$ of pure epimorphisms. Any $f:G\to \colim_i Y_i$ with $G$
finitely presentable lands in some $Y_i$, and then lifts through $X_i$, 
to give a lifting of $f$ itself through $\colim_i p_i$.
\endproof

It now follows, for any class \F, that the  objects injective \wrt \F are closed under products,
retracts, and finite powers. Furthermore, they are closed under $\lambda$-filtered colimits for any regular cardinal $\lambda$ large enough that all domains and codomains of maps in \F are $\lambda$-presentable. Since \V is locally presentable,
such $\lambda$ will exist if \F is small.

It is also convenient to state

\begin{proposition}
For $f:A\to B$ in \K and $C\in\K$, the following are equivalent:
\begin{enumerate}[$(i)$]
\item $C$ is injective \wrt $f$ in \K
\item $G\ct C$ is injective \wrt $f$ in $\K_0$, for all finitely presentable $G$
\item $C$ is injective \wrt $G\cdot f$ in $\K_0$, for all finitely presentable $G$.
\end{enumerate}
\end{proposition}

\noindent
Now turn to weak left adjoints.

\begin{proposition}
Let \A and \K be \V-categories with finite powers, and let $U:\A\to\K$ 
be a \V-functor which preserves finite powers. Then $\eta:K\to UFK$
exhibits $FK$ as a weak left adjoint to $U$ at $K$ if and only if 
it exhibits $FK$ as a weak left adjoint to $U_0:\A_0\to\K_0$ at $K$.
\end{proposition}

\proof
Observe that 
$$\xymatrix{ \A(FK,A) \ar[r]^-{U} & \K(UFK,UA) \ar[rr]^-{\K(\eta,UA)} && \K(K,UA) }$$
is in \E if and only if 
$$\xymatrix{
\V_0(G,\A(FK,A)) \ar[r]^-{\V_0(G,U)} & \V_0(G,\K(UFK,UA)) \ar[rr]^-{\V_0(G,\K(\eta,UA))} &&
\V_0(G,\K(K,UA)) }$$
is surjective, which in turn is the case if and only if 
$$\xymatrix{
\A_0(FK,G\ct A) \ar[r]^-{U_0} & \K_0(UFK,U(G\ct A)) \ar[r]^{\cong} & 
\K_0(UFK,G\ct UA) \ar[rr]^{\K_0(\eta,G\ct UA)} && \K_0(K,G\ct UA) }$$
is surjective.
\endproof

It is also useful to note 

\begin{proposition}
Let $U:\A\to\B$ be a fully faithful functor whose image is closed
under retracts, and which has a weak left adjoint
$\eta B:B\to UFB$ at every object $B\in\B$. Then $C$ is in the image of $U$
if and only if it is injective with respect to all $\eta B$.
\end{proposition}

\proof
If $C$ is in the image of $U$, then injectivity with respect to the $\eta B$
is what it means for the $\eta B$ to give a weak left adjoint. 

Suppose conversely that $C$ is injective with respect to all the $\eta B$.
In particular, it is injective with respect to $\eta C$, and so
$$\xymatrix @C3pc {
\K(UFC,C) \ar[r]^-{\K(\eta C,C)} & \K(C,C) }$$
is in \E. Now $I$ is finitely presentable in $\V_0$, so the identity
$j:I\to\K(C,C)$ lifts to give a map $k:I\to\K(UFC,C)$ which is a retraction
of $\eta C$. This shows that $C$ is a retract of $UFC$ and so is in the
image of $U$.
\endproof

\begin{corollary}
Any weakly reflective subcategory closed under retracts is closed under
products and finite \cotensors.
\end{corollary}
 




Theorem-Schema~\ref{schema:main}, 
in the current setting of \E the pure epimorphisms,
contains \cite[Theorem~4.8]{AR} as the special
case where $\V=\Set$, and indeed our proof of the remaining implication
$(ii)\Rightarrow(iii)$ amounts to reducing 
the general case to the special one, using the existence of
finite \cotensors:

\begin{theorem}
Let \K be a locally presentable \V-category, and \A a full subcategory. 
The following conditions are equivalent:
\begin{enumerate}[$(i)$]
\item \A is the category of objects injective to a small class of maps in \K;
\item \A is accessible, accessibly embedded, and closed under products
and finite \cotensors;
\item \A is accessibly embedded and weakly reflective.
\end{enumerate}
\end{theorem}

\proof
It remains only to prove that $(ii)\Rightarrow(iii)$. Suppose then that 
\A is accessible, accessibly embedded, and closed under products and finite \cotensors.
We must show that it is weakly reflective.

Since $\A_0$ is accessible, and accessibly embedded and closed under products 
in $\K_0$ it follows by \cite[Theorem~4.8]{AR} 
that $\A_0$ is weakly reflective in $\K_0$. Thus for each $K\in\K$
there is an object $K^*\in\A$ and a morphism $r:K\to K^*$ such that for any $A\in\A$,
the morphism $\K_0(r,A)$ is surjective. But \A is closed under finite \cotensors, so 
for any finitely presentable $G\in\V$, the \cotensor $G\ct A$ in \K lies in \A. Thus 
also $\K_0(r,G\ct A)$ is surjective. But $\K_0(r,G\ct A)$ is just $\V_0(G,\K(r,A))$,
and surjectivity of this says that $\K(r,A)$ is a pure epi. This proves that $r:K\to K^*$
is not just a weak reflection of $K$ into $\A_0$ but also an \E-weak reflection of $K$
into \A.
\endproof

Theorem-Schema~B in this setting generalizes \cite[Theorem~4.11]{AR}, and 
follows immediately from the previous theorem and the results of 
Section~\ref{sect:results}:

\enlargethispage{\baselineskip}

\begin{theorem}
For a \V-category \A, the following are equivalent:
\begin{enumerate}[$(i)$]
\item \A 
  is a weakly reflective, accessibly embedded subcategory of $[\C,\V]$ for some small \V-category \C;
\item \A is equivalent to  a small-injectivity class in some locally presentable \V-category;
\item \A is accessible and has products and finite \cotensors;
\item \A is accessible and weakly cocomplete. \endproof
\end{enumerate}
\end{theorem}

\noindent
A \V-category satisfying these conditions is called 
weakly locally presentable, or \E-weakly locally presentable for emphasis.

Theorem-Schema~C generalizes \cite[Theorem~4.13]{AR}:

\begin{theorem}
A \V-category is \E-weakly locally presentable  if and only if it is equivalent to the
\V-category of models of a (limit,\E)-sketch. \endproof
\end{theorem}

We end this section by spelling out a little what \E-weak colimits are in this 
context. Let $S:\D\to\A$ be a \V-functor, and $F:\D\op\to\V$ a presheaf.
An \E-weak colimit of $S$ weighted by $F$ consists of an object
$F*_w S$ and a \V-natural transformation $\eta:F\to\A(S,F*_w S)$ for which
the induced map 
$$\xymatrix{
\A(F*_w S,A)\ar[r] & [\D\op,\V](F,\A(S,A)) }$$
in \V lies in \E for all $A\in\A$. This in turn says that for any finitely presentable
object $G\in\V$ and any $y$ as in the solid part of the diagram 
$$\xymatrix{
\A(F*_w S,A) \ar[r] & [\D\op,\V](F,\A(S,A)) \\
& G \ar[u]_{y} \ar@{.>}[ul]^{x} }$$
there exists a lifting $x$. If \A has finite \cotensors, then this says that any 
\V-natural $F\to\A(S,A)$ arises from some map $F*_w S\to A$ in \A. 

In the case of weak conical colimits, where \D is an ordinary category and
$S:\D\to\A$ an ordinary functor, if \A has finite \cotensors, a weak colimit of 
$S$ in the enriched sense is just a weak colimit in the ordinary sense.

\section{Case 3: \E is the retract equivalences}\label{sect:equivalence}

Here we treat the case $\V=\Cat$, with \E the retract equivalences: these are 
the functors $f:A\to B$ for which there exists a functor $g:B\to A$ with
$fg=1$ and $gf\cong 1$. The fact that the unit object $1$ is \E-projective 
amounts to the fact that retract equivalences are surjective on objects; the 
fact the retract equivalences are cofibrantly generated is well-known: see 
\cite{qm2cat} for example. Thus once again we shall only have to check
the implication $(ii)\Rightarrow(iii)$.

Note that every retract equivalence is in particular a retraction, and so is 
certainly a pure epimorphism. Thus the notion of weakness considered in this
section is ``less weak'' than the notion of weakness for 2-categories arising 
from the pure epimorphisms.

The retract equivalences
are closed under filtered colimits, and they are closed under 
products, \cotensors, and retracts; more generally, they are closed under
{\em flexible limits}, in the sense of \cite{BKPS}. These flexible limits were
shown in \cite{BKPS} to be be all those limits which can be constructed using
products, splitting of idempotents, and two 2-categorical limits called inserters
and equifiers. We have already observed that the retract equivalences are closed 
under products and splittings of idempotents, and it is not too hard to check that
they are also closed under inserters and equifiers. 

There is also another perspective on this, based on the theory developed
in \cite{hty2mnd}. Recall that \Cat has a model structure for which the trivial
fibrations are the retract equivalences and the weak equivalences are the 
equivalences. For any small 2-category \A, the functor 
2-category $[\A\op,\Cat]$ has a ``projective'' model structure for which the trivial
fibrations and weak equivalences are defined ``pointwise'': thus a 2-natural
$\alpha:F\to G$ is a trivial fibration if and only if each component $\alpha A:FA\to GA$
is a trivial fibration in \Cat (that is, a retract equivalence). Now a 2-functor
$E:\A\op\to\Cat$ is cofibrant in this model structure if and only if it is flexible as
a weight, and the fact that the hom 2-functor $[\A\op,\Cat](E,-)$ sends trivial
fibrations in $[\A\op,\Cat]$ to trivial fibrations in \Cat, for cofibrant/flexible $E$, 
is part of the fact that 
the model structure on $[\A\op,\Cat]$ is not just a model category structure
but a model 2-category structure: see \cite{hty2mnd}. Finally to say that
$[\A\op,\Cat](E,-)$ sends pointwise trivial fibrations to trivial fibrations, for all 
cofibrant $E$, is precisely to say that the retract equivalences are closed under
flexible limits.

A related notion is that of pseudolimit. If $F:\D\op\to\Cat$ and $S:\D\to\K$ are
2-functors, the {\em pseudolimit} of $S$ weighted by $F$ is an object 
$\{F,S\}_{ps}$ of \K equipped with a 2-natural isomorphism
$$\K(A,\{F,S\}_{ps})\cong\Ps(\D\op,\Cat)(F\K(A,S))$$
where we have replaced the usual presheaf 2-category $[\D\op,\Cat]$ appearing
in the definition of limit with the 2-category $\Ps(\D\op,\Cat)$ whose objects
are still 2-functors from $\D\op\to\Cat$ but whose morphisms are pseudonatural
transformations, and whose 2-cells are modifications. In general the pseudolimit
$\{F,S\}_{ps}$ is different (non-equivalent) to $\{F,S\}$, but it turns out that 
the pseudolimit $\{F,S\}_{ps}$ can be calculated as an actual weighted limit 
$\{F',S\}$ for a different weight $F'$ (see \cite{BKPS}) and this weight $F'$ is
flexible. Thus pseudolimits are special case of flexible limits. 

Finally, there is the weighted bilimit $\{F,S\}_b$, which is defined by a 
pseudonatural equivalence
$$\K(A,\{F,S\}_b)\simeq\Ps(\D\op,\Cat)(F\K(A,S)).$$
Since every 2-natural isomorphism is a pseudonatural equivalence, a pseudolimit,
if it exists, is also a bilimit. Putting all this together, we see that if a 2-category
has flexible limits, then it has pseudolimits, and so bilimits: see \cite{BKPS} once again.

Turning to our weak notions, we first
consider injectivity. Let $f:A\to B$ be a morphism in 
a \V-category \K. To say that $C$ is \E-injective \wrt $f$, where \E consists of the 
retract equivalences, is
to say that $\K(f,C):\K(B,C)\to\K(A,C)$ is a retract equivalence of categories. 
More explicitly, this means that for each morphism $a:A\to C$ there exists
a $b:B\to C$ with $bf=a$, and for any two $b,b':B\to C$ and any $\alpha: bf\to b'f$
there exists a unique $\beta:b\to b'$ with $\beta f=\alpha$.

Next we turn to a ``weak'' version of the fact that any accessible category 
with limits has an initial object (of course it is in fact cocomplete).
 We shall only need it in the case where the accessible 2-category has flexible 
limits, but it is no harder to prove under the weaker assumption of bilimits.

\begin{lemma}\label{lemma:bi-initial}
Let \A be an accessible 2-category with bilimits. Then \A has a bi-initial
object.
\end{lemma}

\proof 
First we construct the object. 
Let $\lambda$ be a regular cardinal for which \A is $\lambda$-accessible,
and $\A_\lambda$ the full subcategory of $\lambda$-presentable objects. 
Then the bilimit $L$ of the inclusion $J:\A_\lambda\to\A$ will be our 
bi-initial object.

We must show that $\A(L,A)$ is equivalent to the terminal category for all objects
$A$. To see that $\A(L,A)$ is non-empty, observe that any
$A\in\A$ is a $\lambda$-filtered colimit of objects in $\A_\lambda$, so 
in particular, we can find an object $C\in\A_\lambda$ with a morphism 
$f:C\to A$. Then $L$ has a projection $\pi_C:L\to C$ and so we have a
map $f\pi_C:L\to A$. 

{\em Claim: For any $c:C\to L$ with $C$ $\lambda$-presentable, 
$c\pi_C\cong 1_L$.} To see this, observe that if $D$ is any other $\lambda$-presentable
object, then pseudonaturality of the projections gives an isomorphism
$\pi^{\pi_D c}:\pi_D c\pi_C\cong \pi_D=\pi_D 1_L$. These $\pi^{\pi_D c}$ are natural in $D$,
so there is a unique invertible 2-cell $\gamma:c\pi_C\cong 1_L$ with
$\pi_D\gamma=\pi^{\pi_D c}$ for all $D$. {\em This proves the claim.}

We now show that any two objects of $\A(L,A)$ are isomorphic.
Let $g_1,g_2:L\to A$ be any two maps. Since $L$ is also a $\lambda$-filtered
colimit of $\lambda$-presentables, we can find a $\lambda$-presentable 
object $C$ with a morphism $c:C\to L$. Now $g_1 c,g_2 c:C\to A$ are a 
pair of morphisms with $\lambda$-presentable domain, so we can find
a $\lambda$-presentable $D$ with a morphism $d:D\to A$ and 
factorizations $g_1 c=dh_1$ and $g_2 c=dh_2$. Let $\pi_C:L\to C$ 
be the projection of the limit $L$, and observe that by pseudonaturality of the projections
$h_1\pi_C\cong\pi_D\cong h_2\pi_C$; thus $g_1c\pi_C= dh_1\pi_C\cong dh_2\pi_C= g_2 c\pi_C$, and so finally $g_1\cong g_1 c\pi_C\cong g_2 c\pi_C\cong g_2$.

Finally we show that for any two maps $g_1,g_2:L\to A$ there is a unique
2-cell between them. But we already know that all maps $L\to A$ are isomorphic,
so we may as well suppose that both $g_1$ and $g_2$ are given by 
$c\pi_C$ where $c:C\to A$ is some map into $A$ with $\lambda$-presentable
domain. Suppose then that $\phi_1,\phi_2:c\pi_C\to c\pi_C$ are any two 2-cells:
 we must show that $\phi_1=\phi_2$. Let $d:D\to L$ be any map into $L$ with 
$\lambda$-presentable domain. 
Now since $D$ is $\lambda$-presentable, the map $c\pi_C d:D\to A$ and 
2-cells $\phi_1 d,\phi_2 d:c\pi_C d\to c\pi_C d$ factorize through some $e:E\to A$ 
with $E$ a $\lambda$-presentable object,
say as $c\pi_C d=ef$ and $\phi_1 d=e\psi_1$ and $\phi_2 d=e\psi_2$, with
$f:D\to E$ and $\psi_1,\psi_2:f\to f$.

Now $E$, $D$, $f$, $\psi_1$, and $\psi_2$ are all in $\A_\lambda$, and so 
by pseudonaturality of the cone $\pi$, we have $\psi_1\pi_D=\psi_2\pi_D$,
and so $\phi_1 d\pi_D=e\psi_1\pi_D=e\psi_2\pi_D=\phi_2 d\pi_D$;
finally $d\pi_D\cong1$ by the Claim, and so $\phi_1=\phi_2$ as required.
\endproof

\begin{theorem}\label{thm:A-for-Cat}
Let \K be a locally finitely presentable 2-category, and \A a full sub-2-category 
of \K. Then the following conditions are equivalent:
\begin{enumerate}[$(i)$]
\item \A is the category of objects \E-injective to a small class of morphisms in \K;
\item  \A is accessible, accessibly embedded, and closed under flexible limits;
\item \A is accessibly embedded and  \E-weakly reflective.
\end{enumerate}
\end{theorem}

\proof
Note that we have stated $(ii)$ in terms of flexible limits rather than \E-stable 
limits; as observed above, all flexible limits are \E-stable. It remains only to 
prove the implication $(ii)\Rightarrow(iii)$.


Suppose then that \A is accessible, accessibly embedded, and closed under 
flexible limits.  Let $K\in\K$
be given.  Consider the slice 2-category $K/\A$ whose objects are morphisms
$K\to A$ in \K with $A\in\A$. A morphism in $K/\A$ from $f:K\to A$ to 
$g:K\to B$ is a morphism $x:A\to B$ with $xf=g$. A 2-cell from $x$ to $y$ in 
$K/\A$ is a 2-cell $x\to y$ in \K whose restriction along $f$ is the identity.
Now $K/\A$ will have any colimits that \A does, in particular, it will have 
$\lambda$-filtered colimits for any sufficiently large $\lambda$. Similarly,
$K/\A$ has \cotensors since \A does (\cotensors are flexible). Thus 
$K/\A$ will be accessible provided that its underlying ordinary category 
$(K/\A)_0$ is so; but this $(K/\A)_0$ is just the slice category $K/\A_0$ of 
$\A_0$, which is accessible since \A is.

Furthermore, $K/\A$ has flexible limits, since \A and \K do and the inclusion 
preserves them.  It follows by Lemma~\ref{lemma:bi-initial} that $K/\A$ has
a bi-initial object $r:K\to K^*$. The universal property of the bi-initial property 
is that for any object $a:K\to A$ of $K/\A$, the hom-category $(K/\A)((K^*,r),(A,a))$
is equivalent to the terminal category 1. Now this hom-category can be constructed
as a pullback as in 
$$\xymatrix{
(K/\A)((K^*,r),(A,a)) \ar[r] \ar[d] & \K(K^*,A) \ar[d]^{\K(r,A)} \\
1 \ar[r]_{a} & \K(K,A) }$$
Thus the universal property says that the left vertical is an equivalence (and so 
a retract equivalence) and we are to prove that the right vertical is a retract equivalence.

Now retract equivalences are stable under pullback, so if we knew that the right 
vertical were a retract equivalence it would follow immediately that the left vertical 
was one, but here we need to go in the other direction. To do this, we use the 
fact that \A is closed in \K under \cotensors. For any category $G$, the 
functor $[G,\K(r,A)]:[G,\K(K^*,A)]\to[G,\K(K,A)]$ is isomorphic to 
$\K(r,G\ct A):\K(K^*,G\ct A)\to\K(K,G\ct A)$, and any pullback of this 
along a functor $1\to\K(K^*,G\ct A)$ is a retract equivalence. Thus 
$\K(r,A)$ satisfies the conditions of the following lemma, and so is a 
retract equivalence.
\endproof

\begin{lemma}
Let $p:E\to B$ a functor with the property that for every category $C$ 
and every functor $g:C\to B$, if we form the pullback
$$\xymatrix{
P \ar[d]_{q} \ar[r] & [C,E] \ar[d]^{[C,p]} \\
1 \ar[r]_{g} & [C,B] }$$
the resulting functor $q$ is a retract equivalence. Then $p$ is a retract equivalence.
\end{lemma}

\proof
Taking $C=1$ gives the fact that $p$ is surjective on objects, and that if
$e,e'\in E$ with $pe=pe'$ then there is a unique isomorphism $\epsilon:e\cong e'$
sent by $p$ to the identity. 

Suppose  that $e_1,e_2\in E$ with $\beta:pe_1\to pe_2$. Then $\beta$ determines
a map $g:1\to[\two,B]$, and  we can find $\alpha:e'_1\to e'_2$ with $pe'_1=pe_1$,
$pe'_2=pe_2$, and $p\alpha=\beta$. There are now unique isomorphisms 
$\epsilon_1:e_1\cong e'_1$ and $\epsilon_2:e_2\cong e'_2$ sent by $p$ to 
identities, and now the composite
$$\xymatrix{
e_1 \ar[r]^{\epsilon_1} & e'_1 \ar[r]^{\alpha} & e'_2 \ar[r]^{\epsilon^{-1}_2} & e_2 }$$
is sent by $p$ to $\beta$. This proves that $p$ is full.

It remains to show that $p$ is faithful. Suppose then that $\gamma:e_1\to e_2$
is {\em any} morphism in $E$ with $p\gamma=\beta$. We must show
that $\gamma$ equals the composite displayed above, or equivalently that 
$\epsilon_2\gamma=\alpha\epsilon_1$. 

Now $\epsilon_2\gamma$ and $\alpha\epsilon_1$ can be seen as objects of
$[\two,E]$ which are sent by $[\two,p]$ to the same object of $[\two,B]$. Thus
they must be isomorphic, via unique isomorphisms in $[\two,E]$ sent by 
$[\two,p]$ to identities. But such an isomorphism in $[\two,E]$ would 
have components $\theta:e_1\cong e_1$ and $\phi:e'_2\cong e'_2$ satisfying
$\phi\epsilon_2\gamma=\alpha\epsilon_1\theta$ and being sent by $p$
to identities. But then $\theta$ and the identity on $e_1$ are both isomorphisms
$e_1\to e_1$ lying over the identity, so are equal. Similarly $\phi$ is equal to the 
identity,
and it follows that $\epsilon_2\gamma=\alpha\epsilon_1$ as required. 
\endproof

\begin{theorem}
For any 2-category \A, the following conditions are equivalent:
\begin{enumerate}[$(i)$]
\item \A is a 
 weakly reflective, accessibly embedded, full subcategory of some
presheaf 2-category $[\C,\Cat]$;
\item \A is (equivalent to)  a small-injectivity class in some locally finitely presentable 2-category \K;
\item \A is an  accessible 2-category with  flexible limits;
\item \A is  an accessible 2-category with weak colimits.
\end{enumerate}
\end{theorem}


\proof
We have only departed slightly from Theorem-Schema~B, by using flexible
limits rather than \E-stable ones. But this is consistent with the formulation
of Theorem-Schema~A in Theorem~\ref{thm:A-for-Cat} above, and so the 
result follows.




\endproof

A 2-category \A satisfying the conditions of the theorem is called
{\em weakly locally presentable}.

We define a (limit,\E)-sketch to be a small 2-category \C with finite limits,
equipped with a class \F of morphisms. A model of the sketch is a 
finite-limit-preserving 2-functor from \C to \Cat which sends the morphisms
in \F to retract equivalences. (It is also possible, in the usual way, to consider more 
general presentations for such sketches, where we do not assume the existence
of all finite limits.) The models of the sketch are taken to be a full subcategory
of the functor 2-category $[\C,\Cat]$.

\enlargethispage{\baselineskip}
\enlargethispage{\baselineskip}

\begin{remark}
Being a retract equivalence is a purely equational structure: to say that $f$ is
a retract equivalence is to say that there is a $g$ with $fg=1$ and $gf\cong 1$,
and then any 2-functor will send $f$ to a retract equivalence. Thus it might seem
that the class \F makes no difference when it comes to sketching structures.
But there is a subtlety here: if $g$ and the isomorphism $gf\cong 1$ 
were included in the sketch then morphisms would have to be 
strictly natural with respect to them; by merely requiring $f$ to be an equivalence
we only require our morphisms to be strictly natural with respect to $f$. (It
will then follow that they are {\em pseudonatural} with respect to $g$, but 
not necessarily strictly natural.)
\end{remark}


The following theorem follows immediately from the others, as in 
Section~\ref{sect:results}.

\begin{theorem}
A \V-category \A is the \V-category of models of a (limit,\E)-sketch 
if and only if it is a small-injectivity class in a locally presentable 
\V-category \K; in other words, if and only if it is a weakly locally presentable
\V-category. \endproof
\end{theorem}

\section{Examples of weakly locally presentable 2-categories}
\label{sect:Segal}

In this section we focus on the case where $\V=\Cat$ and \E is the class of 
retract equivalences, and exhibit some of the sorts of examples which can arise
as weakly locally presentable 2-categories.

\subsection{2-categories of fibrant objects}

Let \K be a locally presentable 2-category, equipped with a model 2-category 
structure \cite{hty2mnd}; that is, a \Cat-model structure in the sense of \cite{Hovey-book}
for the ``categorical'' or ``natural'' model structure on \Cat. 

Explicitly, this means that there is a model structure on the underlying ordinary
category $\K_0$ satisfying the following condition. Let
$i:A\to B$ be a cofibration and $p:C\to D$ a fibration, and form the pullback
of $\K(i,D):\K(B,D)\to\K(A,D)$ and $\K(A,p):\K(A,C)\to\K(A,D)$. Then the 
induced functor 
$$\xymatrix{
\K(B,C) \ar[r] & \K(A,C)\x_{\K(A,D)} \K(B,D) }$$
is a fibration in \Cat, trivial if either $i$ or $p$ is a weak equivalence.

Now consider the full subcategory \A of \K consisting of the fibrant objects. These 
are the objects $C$ for which $C\to 1$ is a fibration; equivalently, they are 
characterized by the property that for each trivial cofibration $i:A\to B$, the function
$\K_0(i,C):\K_0(B,C)\to\K_0(A,C)$ is surjective, or in other words the functor
$\K(i,C):\K(B,C)\to\K(A,C)$ is surjective on objects. But by the model 2-category 
condition above, this functor is not just surjective on objects but a retract
equivalence.

Thus \A is an \E-injectivity class in \K. If moreover the model structure on \K is
cofibrantly generated, then \A is a small-\E-injectivity class, and so is weakly
locally presentable. 

\subsection{Coflexible presheaves}

For a small 2-category \C, we write $[\C\op,\Cat]$ for the 2-category of 2-functors,
2-natural transformations, and modifications; and we write $\Ps(\C\op,\Cat)$ for
the 2-category of 2-functors, pseudonatural transformations, and modifications.
The inclusion $J:[\C\op,\Cat]\to\Ps(\C\op,\Cat)$ has a left adjoint, sending a
presheaf $F:\C\op\to\Cat$ to a presheaf $F'$ with the property that pseudonatural
transformations from $F$ to $G$ are in bijection with 2-natural transformations
from $F'$ to $G$ (as well as a 2-dimensional aspect to this universal property,
involving the modifications). There is a canonical 2-natural transformation
$q:F'\to F$, which is the component at $F$ of the counit of the adjunction. 
$F$ is flexible when this $q$ has a 2-natural section $s$. It is then a 
consequence that $q$ is a retract equivalence; see \cite{BKPS} for more details,
or \cite{hty2mnd} for the fact that these flexible presheaves are the cofibrant objects 
for the {\em projective}  model structure on $[\C\op,\Cat]$: this is the model
structure for which a 2-natural transformation $f:F\to G$ is a weak equivalence
or a fibration if and only if $fC:FC\to GC$ is one for each object $C$ of \C.

There is also a dual version of these results, using the {\em injective} model
structure \cite[Proposition~A.3.3.2]{HTT} on $[\C\op,\Cat]$ for which $f:F\to G$ is a 
weak equivalence
or a {\em co}fibration if and only if $fC:FC\to GC$ is one for each object $C$ of \C.
The inclusion $J:[\C\op,\Cat]\to\Ps(\C\op,\Cat)$ has a right adjoint, whose image at a presheaf $F$ we shall call 
$\check{F}$, and the component at $F$ of the unit is a 2-natural transformation
$j_F:F\to\check{F}$.  The components $j_FA:FA\to\check{F}A$ of $j_F$ are all
trivial cofibrations in \Cat, and so $j_F$ is a trivial cofibration in $[C\op,\Cat]$ and
an equivalence in $\Ps(\C\op,\Cat)$.
A presheaf $F$ for which this $j_F$ has a 2-natural retraction will
be called {\em coflexible}.

\begin{proposition}
Let \C be a small 2-category.
For a presheaf $F:\C\op\to\Cat$ the following conditions are equivalent:
\begin{enumerate}[(i)]
\item $F$ is coflexible
\item $[\C\op,\Cat](j_F,F):[\C\op,\Cat](\check{F},F)\to[\C\op,\Cat](F,F)$ is 
a retract equivalence
\item $[\C\op,\Cat](j_G,F):[\C\op,\Cat](\check{G},F)\to[\C\op,\Cat](G,F)$ is 
a retract equivalence for all $G$
\item $F$ is fibrant in the injective model structure on $[\C\op,\Cat]$.
\end{enumerate}
\end{proposition}

\proof Here $(iv)\Rightarrow(iii)$ since $j_G$ is a trivial cofibration, and
$(iii)\Rightarrow(ii)$ is trivial. To see that $(ii)\Rightarrow(i)$, observe
that if $[\C\op,\Cat](j_F,F)$ is a retract equivalence then in particular it is surjective
on objects, and so there is some 2-natural $r:\check{F}\to F$ with 
$[\C\op,\Cat](j_F,F)(r)=1$; that is, $rj_F=1$. Thus $F$ is coflexible.

So it remains only to prove that $(i)\Rightarrow(iv)$. Let $F$ be coflexible and 
$j:G\to H$ a trivial cofibration. We must show that for each $u:G\to F$ in $[\C\op,\Cat]$
there exists a $v:H\to F$ in $[\C\op,\Cat]$ with $vj=u$. Since $j$ is a trivial cofibration, there
is a pseudonatural $p:H\to G$ with $pj=1$, and now $up:H\to F$ is pseudonatural
and satisfies $upj=u$. The problem is to replace $up$ by a 2-natural $v:H\to F$
with $vj=u$. 

Let $p_F:\check{F}\to F$ be the pseudonatural map which is the component at $F$
of the counit of the adjunction between $[\C\op,\Cat]$ and $\Ps(\C\op,\Cat)$. 
Then $up$ factorizes as $p_Fv'$ for a unique 2-natural $v':H\to\check{F}$, and 
now $v'j$ and $j_F u$ are 2-natural maps with $p_F v'j=upj=u=p_F j_F u$, and so
$v'j=j_F u$. Finally, since $F$ is coflexible, there is a 2-natural $r:\check{F}\to F$
with $rj_F=1$. Thus $rv'j=rj_F u=u$, and so we may take $v=rv'$.
%
\endproof

Since the injective model structure on $[\C\op,\Cat]$ is cofibrantly generated
(\cite[Appendix~A.3.3]{HTT} again) we are in the situation of the previous section,
and the 2-category $[\C\op,\Cat]\fib$ of coflexible presheaves is weakly locally presentable.

For a more general and more detailed study of coflexibility see the thesis 
\cite{Bourke-thesis}.
One reason to be interested in coflexible presheaves is the following:

\begin{proposition}
For any small 2-category \C, the composite inclusion 
$$\xymatrix{
[\C\op,\Cat]\fib \ar[r] & [\C\op,\Cat] \ar[r] & \Ps(\C\op,\Cat) }$$
is a biequivalence of 2-categories.
\end{proposition}

\proof
First of all, both inclusions are locally fully faithful (fully faithful on the hom-categories).
If $F$ and $G$ are presheaves, with $G$ coflexible, any pseudonatural transformation
$F\to G$ is isomorphic to a 2-natural one; in particular this is the case if $F$ and
$G$ are both coflexible. This proves that the composite is essentially surjective on
the hom-categories.

Finally, in any model category, every object is weakly equivalent to a fibrant one; thus
in $[\C\op,\Cat]$ every presheaf is weakly equivalent to coflexible one;
but weakly equivalent presheaves are pseudonaturally equivalent; that is, 
equivalent in $\Ps(\C\op,\Cat)$. This proves that the composite inclusion is
biessentially surjective on objects, and so a biequivalence.
\endproof

\subsection{Bicategories}

The example of the previous section can be further developed by starting not just 
with all presheaves 
$\C\op\to\Cat$, but only those which preserve some class of limits. This allows
various algebraic structures to be described. Then one could, following
\cite{Segal74}, consider those functors which preserve products only up to homotopy.
For example, if \C has finite coproducts, one could consider functors $F:\C\op\to\Cat$
for which the canonical comparisons $F(C+D)\to FC\x FD$ and $F0\to 1$ are 
retract equivalences. These are the presheaves which are 
\E-injective with respect to the maps $\C(-,C)+\C(-,D)\to\C(-,C+D)$ and 
$0\to C(-,0)$. Once again one would also want to restrict to something like
the coflexible presheaves. 

In this section, however, we have chosen to work through a similar but different
example, involving the structure of bicategory. We start with the 2-category 
$[\Delta\op,\Cat]$, of simplicial
objects in \Cat. It was shown in \cite{2-nerve} that there is a full sub-2-category 
of $[\Delta\op,\Cat]$ which can be identified with the 2-category  \nhom of
bicategories, normal homomorphisms of bicategories, and icons. 

The objects of this full sub-2-category were defined by the following four
requirements, in which we write $X$ for a typical simplicial object in \Cat, and
$X_n$ for the category of $n$-simplices. For each $n$, we may form the 
$n$-fold fibre product $X^n_1$ of $n$ copies of $X_1$ over $X_0$; this represents
the ``composable $n$-tuples'', and comes equipped with a canonical map
$X_n\to X^n_1$ often called the Segal map. The four conditions for $X$ 
to be in the subcategory \nhom are:
\begin{enumerate}[$(i)$]
\item The simplicial object $X$ is 3-coskeletal: this means that it is the right 
Kan extension of its restriction to the full subcategory of $\Delta\op$ 
containing the objects $[0]$, $[1]$, $[2]$, $[3]$;
\item The category $X_0$ of 0-simplices is discrete;
\item The maps $c_2:X_2\to(\Cosk_1 X)_2$ and $c_3:X_3\to(\Cosk_1X)_3$ 
are discrete isofibrations (see below);
\item The Segal map $X_n\to X^n_1$ is a retract equivalences for all $n$.
\end{enumerate}
Here condition $(i)$ says that each $X_n$ with $n>3$ is canonically a limit
of $X_3$, $X_2$, $X_1$, and $X_0$; this is a limit condition so imposing this
restriction does not take us outside of the world of locally presentable categories.
Once again, condition $(ii)$ is a limit condition, since it can be seen as saying 
that the canonical map $X_0\to X^\two_0$ is invertible. A functor 
$f:A\to B$ is called a discrete isofibration if for each object $a\in A$
and each isomorphism $\beta:b\cong fa$ in $B$, there is a unique isomorphism
$\alpha:a'\cong a$ lying over $\beta$. Once again this is a limit condition: it
says that the diagram 
$$\xymatrix{
A^\Iso \ar[r]^{\cod} \ar[d]_{f^{\Iso}} & A \ar[d]^f \\
B^\Iso \ar[r]_{\cod} & B}$$
is a pullback, where
$A^\Iso$ is the category of isomorphisms in $A$, and $\cod$ the 
codomain map. Thus the full sub-2-category of $[\Delta\op,\Cat]$ consisting
of the objects satisfying conditions $(i)$, $(ii)$, and $(iii)$ is still locally presentable.
Finally, condition $(iv)$ is an injectivity condition. For example, to say that
the Segal map $X_2\to X_1\x_{X_0} X_1$ is a retract equivalence is to say that $X$
is \E-injective with respect to the map induced by 
$$\xymatrix{
\Delta_0\cdot I \ar[r]^{\delta_0\cdot I} \ar[d]_{\delta_1\cdot I}  & 
\Delta_1\cdot I \ar[d]^{\delta_0\cdot I} \\
\Delta_1\cdot I \ar[r]_{\delta_2\cdot I}  & \Delta_2\cdot I }$$
 from the pushout of the top and left maps into the bottom corner. 

Thus the 2-category \nhom is weakly locally presentable.

\section{Enriched accessibility and pure subobjects}\label{sect:technical}


In this technical section we state and prove the result, Theorem~\ref{thm:purity} below, which completes the proof of Corollary~\ref{cor:one-two}

We recall from \cite[Section~2D]{AR} that for a locally $\lambda$-presentable ordinary
 category \K, a 
morphism $f:A\to B$ is said to be a {\em $\lambda$-pure monomorphism} if for each commutative square
$$\xymatrix{
C \ar[r]^g \ar[d]_{u} & D \ar[d]^{v} \\
A \ar[r]_f & B}$$
in which $C$ and $D$ are $\lambda$-presentable, there exists a morphism
$w:B\to A$ with $wg=u$. These $\lambda$-pure monomorphisms are always
monomorphisms; in \cite{AR} they were called simply $\lambda$-pure morphisms,
but since pure epimorphisms are generally  not pure morphisms in this sense,
we have chosen to use the name $(\lambda)$-pure monomorphism to reduce the
risk of confusion. The $\lambda$-pure monomorphisms
can be characterized as the closure in $\K^\two$ of the 
split monomorphisms under $\lambda$-filtered colimits. 

From the notion of $\lambda$-pure monomorphism we obtain the notion of 
$\lambda$-pure subobject, which 
plays an important role in the theory of accessible categories. 
Here we shall use this notion in the context of enriched categories; the notion 
itself remains unchanged in our enriched context. First we state 
\cite[Theorem~2.34]{AR}, combined with the remark that immediately follows it as:

\begin{theorem}[Ad\'amek-Rosick\'y]
For any $\lambda$-accessible category \K, the (non-full) subcategory 
$\Pure_\lambda\K$ of \K consisting of all objects and all $\lambda$-pure 
monomorphisms
is accessible, and closed in \K under $\lambda$-filtered colimits.   
\end{theorem}

We apply this to the case of a locally $\lambda$-presentable 
\V-category \K: then $\Pure_\lambda\K_0$ is accessible, and closed in $\K_0$
under $\lambda$-filtered colimits. Since not just $\K_0$ but \K has $\lambda$-filtered
colimits, the inclusion $\Pure_\lambda\K_0$ in $\K_0$ sends $\lambda$-filtered
colimits in $\Pure_\lambda\K_0$ to $\lambda$-filtered colimits in~\K.

\begin{theorem}\label{thm:purity}
Let \K be a locally $\lambda$-presentable \V-category and \A a full subcategory closed
under $\lambda$-filtered colimits. If $\A_0$ is closed in $\K_0$ under 
$\lambda$-pure subobjects, then \A is accessible.
\end{theorem}

\proof
The proof follows that of \cite[Corollary~2.36]{AR}, merely taking care to use
enriched notions where necessary. By the previous theorem, we know that 
$\Pure_\lambda\K_0$ is accessible, and that $\lambda$-filtered colimits in
$\Pure_\lambda\K_0$ are $\lambda$-filtered colimits in \K. Let $\mu_0$
be some regular cardinal greater than or equal to $\lambda$ for which 
$\Pure_\lambda\K_0$ is $\mu_0$-accessible. Now choose a regular 
cardinal $\mu\triangleright\mu_0$ such that each $\mu_0$-presentable
object in $\Pure_\lambda\K_0$ is $\mu$-presentable in \K: this is possible
since \K is locally presentable, and the set of all $\mu_0$-presentable objects 
in $\Pure_\lambda\K_0$ is small. Then each $\mu$-presentable object of 
$\Pure_\lambda\K_0$ is a $\mu$-small $\mu_0$-filtered
colimit of $\mu_0$-presentable objects, and so is a $\mu$-small colimit in \K
of $\mu$-presentable objects, and so is $\mu$-presentable in \K.

For any object $A\in\A$, we can write $A$ as a $\mu$-filtered colimit in 
$\Pure_\lambda\K_0$ of $\mu$-presentable objects. This is equally a 
$\mu$-filtered colimit in \K of $\mu$-presentable objects (in \K!) Furthermore,
each vertex of the diagram is a $\lambda$-pure subobject of $A$, so is in \A, 
and so finally we have written $A$ as a $\mu$-filtered colimit in \A of 
$\mu$-presentable objects. Since \A has $\mu$-filtered colimits (and even
$\lambda$-filtered colimits) it follows that \A is $\mu$-accessible.
\endproof


\bibliographystyle{plain}


\end{document}